\documentclass[twoside,leqno,A4,10pt]{amsart}
\usepackage{amsfonts}
\usepackage{amsmath}
\usepackage{amscd}
\usepackage{amssymb}
\usepackage{amsthm}
\usepackage{amsrefs}
\usepackage{latexsym}
\usepackage{bbm}
\usepackage{manfnt}
\numberwithin{equation}{section}

\begin{document}
\newtheorem*{theorem}{Theorem}
\newtheorem{lemma}{Lemma}
\newtheorem*{corollary}{Corollary}
\numberwithin{equation}{section}
\newcommand{\dif}{\mathrm{d}}
\newcommand{\intz}{\mathbb{Z}}
\newcommand{\ratq}{\mathbb{Q}}
\newcommand{\natn}{\mathbb{N}}
\newcommand{\comc}{\mathbb{C}}
\newcommand{\rear}{\mathbb{R}}
\newcommand{\prip}{\mathbb{P}}
\newcommand{\uph}{\mathbb{H}}
\newcommand{\fief}{\mathbb{F}}
\newcommand{\majorarc}{\mathfrak{M}}
\newcommand{\minorarc}{\mathfrak{m}}
\newcommand{\sings}{\mathfrak{S}}

\title{On the sum of integers from some multiplicative sets and some powers of integers}
\author{Timothy Foo}
\date{\today}
\maketitle

\begin{abstract}
We show that if there exists an integer subject to some congruence conditions that cannot be written as the sum of the norm of an ideal in $\mathbb{Z}[\exp(2\pi i/2^k)]$ and at most $k$ powers of $2$, $k\geq 3$, then there are infinitely many integers such integers. Also, if there exists an integer that cannot be written as the sum of an integer which is the norm of an ideal in in $\mathbb{Z}[\exp(2\pi i/p)]$ and at most $p-2$ powers of $p$, where $p\geq 3$ is a prime, then there are infinitely many such integers. Finally it is shown that there are infinitely many integers not the sum of the norm of an ideal in $\mathbb{Z}[\exp(2\pi i/p)]$ and at most $p-2$ powers of $p$, for $p\geq 3$ prime.
\end{abstract}


\section{Introduction}

Crocker [1] proved that if there exists an integer $N_0\equiv 0\bmod 36$ that is not expressible as the sum of two squares and at most two powers of $2$, then there are infinitely many. He proved that there are infinitely many integers not the sum of two squares and at most two powers of $2$. Platt and Trudgian [3] gave a shorter proof and showed that if there exists an integer $N_0\equiv 0\bmod 18$ that is not expressible as the sum of two squares and at most two powers of $2$, then there are infinitely many. They then showed the existence of such an $N_0$, thus showing that there are infinitely many integers not representable as the sum of two squares and at most two powers of $2$. In this note, we show that if there exists an integer $N_0$ such that $N_0\equiv 0 \bmod 196$ and $N_0$ is not expressible as the sum of $N(x,y,z,w)$ and three powers of $2$, where $N(x,y,z,w)=x^4+(4wy+2z^2)x^2+(-4zy^2+4w^2z)x+(y^4+2w^2y^2-4wz^2y +z^4+w^4)$, $x,y,z,w\in\mathbb{Z}$, then there are infinitely such integers. $N(x,y,z,w)$ is the norm of the element $x+y\zeta + z\zeta^2+w\zeta^3$ in $\mathbb{Z}[\zeta]$, where $\zeta=\exp(2\pi i/8)$. In general, we show that for $k\geq 3$, if there exists an integer $N_0$ such that $N_0\equiv 0\bmod 2^{k-1}(2^k-1)^2$ and $N_0$ cannot be represented as the sum of an integer which is the norm of an ideal in $\mathbb{Z}[\exp(2\pi i/2^k)]$ and at most $k$ powers of $2$, then there are infinitely many such integers. Also, if there exists an $N_0$ such that $N_0$ cannot be expressed as the sum of an integer which is the norm of an ideal in $\mathbb{Z}[\exp(2\pi i/p)]$ and at most $p-2$ powers of $p$, where $p\geq 3$ is a prime, then there are infinitely many such integers.
\newline
\begin{lemma} \label{1}
Let  $N(x,y,z,w)=x^4+(4wy+2z^2)x^2+(-4zy^2+4w^2z)x+(y^4+2w^2y^2-4wz^2y +z^4+w^4)$, $x,y,z,w\in\mathbb{Z}$. Suppose that an integer $n$ cannot be represented by $N(x,y,z,w)$. Then neither can $2^an$ for $a\geq 0$.
\end{lemma}
\noindent
{\it Proof.} Since $n$ cannot be represented by $N(x,y,z,w)$, its prime factorisation must contain an odd prime $p\not\equiv 1 \bmod 8$ to an odd power. If the prime factorisation of $n$ contains an odd prime $p\not\equiv 1 \bmod 8$ to an odd power, then so does that of $2^an$, $a\geq0$.
\newline
\begin{lemma} \label{2}
Let  $N(x,y,z,w)=x^4+(4wy+2z^2)x^2+(-4zy^2+4w^2z)x+(y^4+2w^2y^2-4wz^2y +z^4+w^4)$, $x,y,z,w\in\mathbb{Z}$. Suppose that an integer $n\equiv 0\bmod 196$ cannot be represented as the sum of $N(x,y,z,w)$ and at most three powers of $2$. Then neither can $2^an$ for $a\geq 0$.
\end{lemma}
\noindent
{\it Proof.} By assumption, $n$, $n-2^a$, $n-2^a-2^b$ and $n-2^a-2^b-2^c$, $a,b,c\geq 0$, are not representable by $N(x,y,z,w)$, $x,y,z,w\in\mathbb{Z}$. Therefore, by Lemma 1, $2n$, $2n-2^a$, $2n-2^a-2^b$ and $2n-2^a-2^b-2^c$, $a,b,c\geq 1$, are not representable by $N(x,y,z,w)$, $x,y,z,w\in\mathbb{Z}$. Also, since $n\equiv 0\bmod 4$, $2n-1$, $2n-1-2^a$ with $a\geq1$, $2n-1-1-1$, $2n-1-2-2$, $2n-1-2^a-2^b$ with $a,b\geq 2$, $2n-1-2-2^a$ with $a\geq 3$ are all odd but not $1\bmod 8$, so they are not representable by $N(x,y,z,w)$, $x,y,z,w\in\mathbb{Z}$. On the other hand, the cases $2n-1-1$ and $2n-1-1-2^a$ with $a\geq 1$ reduce to the case with one and two powers of $2$ respectively. Finally, since $n\equiv 0\bmod 49$, $2n-1-2-4=2n-7\equiv 42\bmod 49$ and hence is not representable by $N(x,y,z,w)$, $x,y,z,w\in\mathbb{Z}$.
\newline\newline
\section{Generalisation}
Let $\mathcal{N}_k\subset\{x^2+y^2| x,y,\in\mathbb{Z}\}\subset\mathbb{Z}$ be the set of integers that are norms of ideals in $\mathbb{Z}[\exp(2\pi i/2^k)]$. This set of integers has the property that a power of an odd prime, $p^a$, divides $m\in \mathcal{N}_k$ if and only if $p^a\equiv 1 \bmod 2^k$. Therefore, all odd integers in $\mathcal{N}_k$ are $1\bmod 2^k$. Let $ord(p,2^k)$ denote the order of $p$ in $(\mathbb{Z}/2^k\mathbb{Z})^{*}$. The density of the set, for $k>2$, is
$$
\prod_{d|2^{k-2}}\prod_{ord(p,2^k)=d}\frac{1-1/p}{1-1/p^d}=0.
$$ 
Let $\mathcal{N}_k(x)=|\{n\in\mathcal{N}_k|0<n<x\}|$. By Theorem 1 of [2], $\mathcal{N}_k(x)\asymp \frac{x}{(\log x)^{1-1/2^{k-1}}}$.
\begin{lemma} \label{3}
Let $k\geq 3$. Suppose that an integer $n\not\in \mathcal{N}_k$. Then neither is $2^an$ for $a\geq 0$.
\end{lemma}
\noindent
{\it Proof.}  Since $\mathcal{N}_k$ has the property that a power of an odd prime, $p^a$, divides $m\in \mathcal{N}_k$ if and only if $p^a\equiv 1 \bmod 2^k$, then $n\not\in \mathcal{N}_k$ implies $2n\not\in \mathcal{N}_k$.
\newline\newline
A generalisation of Lemma 2 holds for all $k>3$. 
\begin{lemma}\label{4}
Let $k\geq 3$. Suppose that an integer $n\equiv 0\bmod 2^{k-1}(2^k-1)^2$ cannot be represented as the sum of an integer in $\mathcal{N}_k$ and at most $k$ powers of $2$. Then neither can $2^an$ for $a\geq 0$.
\end{lemma}
\noindent
{\it Proof.} Suppose that $n\equiv 0\bmod 2^{k-1}$ is not representable as the sum of an integer in $\mathcal{N}_k$ and at most $k$ powers of $2$. By inspection, the only case that needs to be checked is whether $2n-(2^k-1)$ is in $\mathcal{N}_k$, since any sum of $k$ powers of $2$ not equal to $2^k-1$ is either even or not $2^k-1\bmod 2^k$. For each $k$, although $2n-(2^k-1)\equiv 1\bmod 2^k$, let $n\equiv 0 \bmod (2^k-1)^2$. Then $2n-(2^k-1)\equiv -(2^k-1) \bmod (2^k-1)^2$. Since $2^k-1\equiv 3 \bmod 4$, there exists some $q^m$, $q\equiv 3 \bmod 4$ prime and $m$ odd, dividing $2^k-1$, and $((2^k-1)/q^m,q)=1$. Then $2n-(2^k-1)\equiv -(2^k-1)=-q^m((2^k-1)/q^m) \bmod q^{2m}$. However, a sum of two squares must be congruent to $q^{m^{\prime}}m^{\prime\prime}\bmod q^{2m}$ for some even $m^{\prime}$ and $(m^{\prime\prime},q)=1$ when $q$ is a prime such that $q\equiv 3 \bmod 4$, therefore $2n-(2^k-1)$ is not the sum of two squares, and is not in $\mathcal{N}_k$.
\newline\newline
Let $\mathcal{M}_p$ be the set of positive integers that are norms of ideals of $\mathbb{Z}[\exp(2\pi i/p)]$, $p\geq 3$ a prime. This set of integers has the property that for a prime $q\not=p$, a power of $q$, $q^a$ divides an integer in $\mathcal{M}_p$ if and only if $q^a\equiv 1 \bmod p$. The integers in $\mathcal{M}_p$ that are relatively prime to $p$ are all $1\bmod p$.\\By Theorem 1 of [2], $|\{n\in \mathcal{M}_p| 0< n< x\}| \asymp \frac{x}{(\log x)^{1-1/(p-1)}}$.
\newline
\begin{lemma}\label{5}
Let $p\geq 3$ be a prime. Suppose that an integer $n\not\in \mathcal{M}_p$. Then neither is $p^an$ for $a\geq 0$.
\end{lemma}
\noindent
{\it Proof.}  Since $\mathcal{M}_p$ has the property that a power of a prime $q\not=p$, $q^a$, divides $m\in \mathcal{M}_p$ if and only if $q^a\equiv 1 \bmod p$, then $n\not\in \mathcal{M}_p$ implies $pn\not\in \mathcal{M}_p$.
\newline
\begin{lemma}\label{6}
Let $p\geq 3$ be a prime. Suppose that an integer $n$ cannot be represented as the sum of an integer in $\mathcal{M}_p$ and at most $p-2$ powers of $p$. Then neither can $p^an$ for $a\geq 0$.
\end{lemma}
\noindent
{\it Proof.} Suppose $n$ meets the conditions of the lemma. Then $pn-($ sum of $p-2$ powers of $p)$ is either $0\bmod p$, where we use Lemma 5, or neither $0$ nor $1\bmod p$.
\newline
\begin{lemma}\label{7}
$11$ is not the sum of an integer in $\mathcal{M}_3$ and at most one power of $3$. $9$ is not the sum of an integer in $\mathcal{M}_5$ and at most three powers of $5$. $20$ is not the sum of an integer in $\mathcal{M}_7$ and at most five powers of $7$.
\end{lemma}
\noindent
By the previous two lemmas, there are infinitely many integers not the sum of the norm of an ideal in $\mathbb{Z}[\exp(2\pi i/p)]$ and at most $p-2$ powers of $p$ for $p=3,5,7$. It is also true that there are infinitely many integers not the sum of the norm of an ideal in $\mathbb{Z}[\exp(2\pi i/p)]$ and at most $p-2$ powers of $p$ for $p\geq 3$ prime.
\begin{theorem}
There are infinitely many integers not the sum of the norm of an ideal in $\mathbb{Z}[\exp(2\pi i/p)]$ and at most $p-2$ powers of $p$, for $p\geq 3$ prime.
\end{theorem}
\noindent
{\it Proof.} Suppose $n\equiv -1\bmod p$. The sum of $p-2$ powers of $p$ is never $-1\bmod p$ and is exactly $p-2$ if it is $-2\bmod p$. Integers in $\mathcal{M}_p$ are either $0$ or $1\bmod p$. Therefore, if $n$ is the sum of an integer in $\mathcal{M}_p$ and $p-2$ powers of $p$, it is the sum of a $1\bmod p$ integer in $\mathcal{M}_p$ and $p-2$. There are infinitely many numbers $q_1q_2\equiv 1 \bmod p$ where $q_1,q_2\not\equiv 1\bmod p$ are primes, and $q_1q_2\not\in\mathcal{M}_p$. In these cases, $n=q_1q_2+p-2$ is not the sum of an integer in $\mathcal{M}_p$ and $p-2$ powers of $p$.

\section{Computational aspects}
The following PARI/GP program may be used to find a small integer not the sum of the norm of an ideal in $\mathbb{Z}[\exp(2\pi i/p)]$ and at most $p-2$ powers of $p$, for $p=7$. It may be modified for other primes.
\newline
\begin{verbatim}
p=7;v=vector(30);for(n=1,30,v[n]=1);
for(n=1,30,a=factor(n);m=matsize(a);
for(i=1,m[1],if(a[i,1]%p!=0%p,if((a[i,1]^a[i,2])%p!=1%p,v[n]--))));
for(n=1,30,if(v[n]<0,v[n]=0));
v2=vector(30);for(i=0,10,if(p^i<30,v2[p^i+1]++));v2[1]=1;
v3=vector(30);for(n=1,30,for(n1=1,30,for(n2=1,30,for(n3=1,30,for(n4=1,30,for(n5=1,30, 
if(n+n1+n2+n3+n4+n5<=30+p-2, 
v3[n+n1+n2+n3+n4+n5-(p-2)]+=v[n]*v2[n1]*v2[n2]*v2[n3]*v2[n4]*v2[n5] )))))))
\end{verbatim}
\vspace*{.7cm}
We are not sure whether or not there exists some $N_0\equiv 0 \bmod 2^{k-1}(2^k-1)^2$ that cannot be represented as the sum of an integer in $\mathcal{N}_k$ and at most $k$ powers of $2$, $k\geq 3$. Since
$$
\left(\sum_{\substack{a<x\\a\in \mathcal{N}_k}}1\right)\left(\sum_{1\leq 2^b<x}1\right)^k\asymp x(\log x)^{k-1+1/2^{k-1}},
$$
we don't know whether or not there are infinitely many integers not representable as the sum of an integer in $\mathcal{N}_k$ and $k$ powers of $2$, $k\geq 3$. Also,
$$
\left(\sum_{\substack{a<x\\a\in \mathcal{M}_p}}1\right)\left(\sum_{1\leq p^b<x}1\right)^{p-2}\asymp x(\log x)^{p-3+1/(p-1)}.
$$
\newline\newline
Platt and Trudgian [3] showed the interesting observation that there are a positive proportion of integers not the sum of two squares and at most one power of $2$, and both Crocker [1] and Platt and Trudgian [3] proved that there are infinitely many integers not the sum of two squares and at most two powers of $2$. The following may also be observed on sums of integers in $\mathcal{N}_k$ and at most $k-1$ powers of $2$, $k\geq 3$ or sums of integers in $\mathcal{N}_k$ and at most $k$ powers of $2$, $k\geq 3$. By Theorem 1 of [2],
$$
|\{n\in\mathcal{N}_k| n<x, n\equiv a\bmod 2^k \}| \mbox{ is } \begin{cases}\ll \frac{x}{(\log x)^{1-1/2^{k-1}}},\mbox{if }a=1 \mbox{ or a power of }2\bmod 2^k\\
\mbox{finite otherwise.}
\end{cases}
$$
Also, the number of integers that are the sum of at most $k-1$ powers of $2$ and less than $x$ and that are congruent to $a\bmod 2^k$ is finite if and only if the binary representation of $a$ has $k-1$ or more $1$s. Therefore, the number of integers less than $x$ and congruent to $2^k-1\bmod 2^k$ that are the sum of an integer from $\mathcal{N}_k$ and at most $k-1$ powers of $2$, is $\ll \frac{x}{(\log x)^{1-1/2^{k-1}}}$. Therefore, a positive proportion of integers less than $x$ and congruent to $2^k-1\bmod 2^k$ are not the sum of an integer from $\mathcal{N}_k$ and at most $k-1$ powers of $2$, $k\geq 3$. When there are at most $k$ powers of $2$, the number of integers that are the sum of at most $k$ powers of $2$ and less than $x$ that are congruent to $a\bmod 2^k$ is $\gg \log x$ if the binary expansion of $a$ has exactly $k-1$ $1$s, and finite if it has $k$ $1$s. The number of pairs $(n_1,n_2)$ with $n_1,n_2$ less than $x$, $n_1\in\mathcal{N}_k$, $n_2$ the sum of at most $k$ powers of $2$ and congruent to some $a\bmod 2^k$ where the binary expansion of $a$ has exactly $k-1$ $1$s, is $\gg x(\log x)^{1/2^{k-1}}$.
\newline\newline
\noindent{\bf Acknowledgements.}  
\newline\newline\noindent
The author thanks Tim Trudgian for helpful conversations.

\section*{References}
\noindent
[1] R. C. Crocker. On the sum of two squares and two powers of $k$. Colloq. Math. 112 (2008) 235--267.
\newline\newline
[2] P. Moree. Counting numbers in multiplicative sets: Landau versus Ramanujan. arXiv:1110.0708.
\newline\newline
[3] D. Platt, T. Trudgian. On the sum of two squares and at most two powers of $2$. arXiv:1610.01672. 

\bibliography{biblio}
\bibliographystyle{amsxport}

\vspace*{.7cm}

\noindent\begin{tabular}{p{8cm}p{8cm}}
Deptartment of Mathematical and Computational Sciences\\
University of Toronto Mississauga \\
Email: {\tt tch.foo@utoronto.ca} \\
\end{tabular}

\end{document}